
\def\mod{\rm mod}
\def\Z{\bf Z}

\baselineskip=14pt
\parskip=10pt

\magnification=\magstephalf

\def\1{{\overline{1}}}
\def\2{{\overline{2}}}
\parindent=0pt
\overfullrule=0in

\def\frac#1#2{{#1 \over #2}}


\bf
\centerline
{
A Meta-Algorithm for Creating Fast Algorithms for Counting ON Cells 
}
\centerline
{
 in Odd-Rule Cellular Automata}
\rm
\bigskip
\centerline
{\it By Shalosh B. EKHAD, N.~J.~A. SLOANE, and  Doron ZEILBERGER}

{\bf Abstract}: By using the methods of Rowland and Zeilberger (2014),
we develop a meta-algorithm that, given a polynomial (in one or more variables), 
and a prime $p$, produces
a fast (logarithmic time) algorithm that takes a positive integer $n$ and 
outputs the number of times each residue class modulo $p$
appears as a coefficient when the polynomial is raised to the power $n$
and the coefficients are read modulo $p$. 
When $p=2$, this has applications to counting the ON cells in certain ``Odd-Rule''
cellular automata.
(This article is accompanied by a Maple package, {\tt CAcount}, 
as well as numerous examples of input and output
files, all of which can be obtained from the web page for this article: \hfill\break
{\tt http://www.math.rutgers.edu/\~{}zeilberg/mamarim/mamarimhtml/CAcount.html}).

{\bf Preface}

The number of ON cells in the $n$th generation of
an ``Odd-Rule'' cellular automaton is found by raising the defining
polynomial (in which the number
of variables is equal to the dimension of the ambient space)
to the $n$th power, reading the coefficients modulo $2$, and counting the remaining monomials---or
equivalently, setting all the variables equal to $1$ (see [Sl]
for a detailed discussion).

The purpose of this article is to describe a meta-algorithm, inspired by
a recent paper of Eric Rowland and Doron Zeilberger [RZ],  that takes such a polynomial
as input,
and outputs a recurrence scheme
 that enables the fast (logarithmic time)  computation of terms of the
sequence giving the number of ON cells at time $n$.
This provides an alternative, computer proof of Theorems
4 and 5 of [Sl].

{\bf A toy example}

Following the {\it Gelfand Principle}, let's illustrate the method with a simple example
that can  be done by hand.
We will later describe how this method can be `taught' to a computer,
which will  then be able to do far more complicated cases, impossible for humans.

Consider the sequence
$$
a_1(n):=(1+x+x^2)^n \, {\mod} \, 2 \, \Bigl \vert_{x =1} \quad,
$$
(sequence A071053 in [OEIS]),
and suppose we wish to compute $a_1(10^{100})$, 
or  $a_1(n)$ for any very large $n$.

Of course, direct computation is hopeless, even if we reduce modulo 2 at
each step  and use the
repeated squaring trick that makes RSA possible 
($P^n=(P^{n/2})^2$ if $n$ is even,  $P^n=PP^{n-1}$ if $n$ is odd),
since the polynomials, before we set $x=1$, are far too big for our modest universe.
What we will do is adapt this trick so that we can also make
the substitution $x=1$ at 
intermediate steps.

First let's try to relate $a_1(2n)$ to $a_1(n)$, using the {\bf Freshman's Dream} identity $P(x)^p \equiv P(x^p) \, {\mod} \, p$:
$$
a_1(2n) =(1+x+x^2)^{2n} \, {\mod} \, 2 \, \Bigl \vert_{x=1} = ((1+x+x^2)^{2})^{n} \, {\mod} \, 2 \, \Bigl \vert_{x=1} 
$$
$$
= (1+x^2+x^4)^{n} \, {\mod} \, 2 \, \Bigl \vert_{x=1}  = (1+x+x^2)^{n} \, {\mod} \, 2 \, \Bigl \vert_{x=1}  
\eqno(EvenCase1)
$$
(replacing $x^2$ by $x$). Hence
$$
a_1(2n)=a_1(n) \quad .
\eqno(Recurrence1even)
$$

Now we do the same thing for $a_1(2n+1)$:
$$
a_1(2n+1)\, =\, (1+x+x^2)^{2n+1} \, {\mod} \, 2 \, \Bigl \vert_{x=1} = (1+x+x^2) \, ((1+x+x^2)^{2})^{n} \, {\mod} \, 2 \, \Bigl \vert_{x=1} 
$$
$$
\, = \, (1+x+x^2) \, (1+x^2+x^4)^{n} \, {\mod} \, 2 \, \Bigl \vert_{x=1} 
$$
$$
\, = \, (1+x^2) \, (1+x^2+x^4)^{n} \, {\mod} \, 2 \, \Bigl \vert_{x=1} \, + \, x \, (1+x^2+x^4)^{n} \, {\mod} \, 2 \, \Bigl \vert_{x=1}  \quad .
\eqno(OddCase1)
$$
In the first term, once again, we can replace $x^2$ by $x$, getting an {\bf uninvited guest},  
$a_2(n)$, say:
$$
a_2(n):=(1+x) \, (1+x+x^2)^n \, {\mod} \, 2 \, \, \Bigl \vert_{x=1} \quad .
$$
As for the second term of Eq. $(OddCase1)$, multiplying by $x$ does not change anything, so this is equal to
$(1+x^2+x^4)^{n} \, {\mod} \, 2 \, \Bigl \vert_{x=1}$, which, again replacing $x^2$ by $x$, 
is our old friend $a_1(n)$.
Hence
$$
a_1(2n+1)=a_2(n)+a_1(n) \quad .
\eqno(Recurrence1odd)
$$
But this pair of recurrences is useless unless we can handle $a_2(n)$. So let's try 
the same technique on it.
A priori, this may force us to introduce terms $a_3(n)$, $a_4(n)$, etc., 
and lead us into an infinite regression,
also known as a {\it Ponzi scheme}, but let's hope for the best.

Again we start with  $a_2(2n)$.  Using the Freshman's Dream, and the fact that
multiplying a polynomial by $x$ (or any other monomial) does not
affect the result if we are going to read it modulo $2$ and set $x=1$, we have
$$
a_2(2n) \, =\, (1+x) \, (1+x+x^2)^{2n} \, {\mod} \, 2 \, \Bigl \vert_{x=1} \, = \,
(1+x) \cdot ((1+x+x^2)^2)^n \, {\mod} \, 2 \, \Bigl \vert_{x=1} 
$$
$$
=\, (1+x) \cdot (1+x^2+x^4)^n \, {\mod} \, 2 \, \Bigl \vert_{x=1} 
\, = \, 1 \cdot (1+x^2+x^4)^n \, {\mod} \, 2 \, \Bigl \vert_{x=1} \, + \, x \, (1+x^2+x^4)^n \, {\mod} \, 2 \, \Bigl \vert_{x=1} \,
$$
$$
=\, 2 \, (1+x^2+x^4)^n \, {\mod} \, 2 \, \Bigl \vert_{x=1} \,
=\, 2 \, (1+x+x^2)^n \, {\mod} \, 2 \, \Bigl \vert_{x=1} \, = \, 2a_1(n) \quad .
$$
Hence
$$
a_2(2n)=2a_1(n) \quad .
\eqno(Recurrence2even)
$$

Now for $a_2(2n+1)$. We have
$$
a_2(2n+1) \, = \, (1+x) \cdot (1+x+x^2)^{2n+1} \, {\mod} \, 2 \, \Bigl \vert_{x=1} 
$$
$$ 
= \,((1+x) \cdot (1+x+x^2)) \cdot ((1+x+x^2)^2)^n \, {\mod} \, 2 \, \Bigl \vert_{x=1} \, 
$$
$$
= \, (1+2x+2x^2+x^3) \cdot (1+x^2+x^4)^n \, {\mod} \, 2 \, \Bigl \vert_{x=1} \, 
$$
$$
= \,
(1+x^3) \cdot (1+x^2+x^4)^n \, {\mod} \, 2 \, \Bigl \vert_{x=1} 
$$
$$
= \, 1 \cdot (1+x^2+x^4)^n \, {\mod} \, 2 \, \Bigl \vert_{x=1} \, + \, x^3 \cdot (1+x^2+x^4)^n \, {\mod} \, 2 \, \Bigl \vert_{x=1}
$$
$$
=\, (1+x^2+x^4)^n \, {\mod} \, 2 \, \Bigl \vert_{x=1} \, + \, (1+x^2+x^4)^n \, {\mod} \, 2 \, \Bigl \vert_{x=1}
$$
$$
= \,
(1+x+x^2)^n \, {\mod} \, 2 \, \Bigl \vert_{x=1} \, + \, (1+x+x^2)^n \, {\mod} \, 2 \, \Bigl \vert_{x=1}= 2a_1(n) \quad .
$$
Hence
$$
a_2(2n+1)=2a_1(n) \quad .
\eqno(Recurrence2odd)
$$
So the {\it uninvited guest}, $a_2(n)$,  did not invite further guests, and now we have a super-fast way to compute
$a_1(n)$ for large $n$, using the {\bf system} 
$$
a_1(2n)=a_1(n) \quad , \quad a_1(2n+1)=a_1(n) +a_2(n) \quad ;
$$
$$
a_2(2n)=2a_1(n) \quad , \quad a_2(2n+1)=2a_1(n) \quad .
\eqno(System)
$$

For certain ``odd-rule'' cellular automata, the sequence $a_1(n), n \ge 0$
is completely determined by the subsequence $b_1(k):=a_1(2^k-1), k \ge 0$ [Sl],
and the $b_1(k)$, unlike the $a_1(n)$, often have simple generating functions,
which we can derive (rigorously) by these methods. 
With $a_1(n)$ as defined above, let
$$
f_1(t):= \sum_{k=0}^{\infty} b_1(k) t^k
$$
be the generating function for $b_1(k)$, and similarly define
 $b_2(k):=a_2(2^k-1)$ and
$$
f_2(t):= \sum_{k=0}^{\infty} b_2(k) t^k \quad .
$$
From  Eq. $(System)$, we have
$$
b_1(k)=b_1(k-1)+ b_2(k-1) \quad , \quad b_2(k)= 2b_1(k-1) \quad ,
$$
and since by direct computation, $b_1(0)=1$, $b_2(0)=2$, 
we arrive at a system of two linear equations for the {\it unknowns} $f_1(t)$ and $f_2(t)$:
$$
\{ \, f_1(t)=1+tf_1(t) + t f_2(t) \quad, \quad f_2(t)=2+ 2t f_1(t) \, \} \quad ,
$$
whose solution is
$$
f_1(t)= \frac{1 \, + \, 2t}{(1+t)(1-2t)} \quad , \quad
f_2(t)= \frac{2}{(1+t)(1-2t)} \quad 
$$
(A001045, A014113 in [OEIS]). But we really don't care about $f_2(t)$, we just needed it in order to find $f_1(t)$, so {\it now} we can safely discard it,
and get the 

{\bf Theorem:} 
$$
f_1(t)= \frac{1 \, + \, 2t}{(1+t)(1-2t)} \quad .
$$

{\bf The  general case}

Fix once and for all a prime $p$ and a polynomial $P=P(x_1, \dots, x_k) \in {\Z}[x_1,\ldots,x_k]$. 
If $A(x_1, \dots, x_k)$ is any element of ${\Z}[x_1,\ldots,x_k]$, we define the 
 {\it functional} 
$$
A(x_1, \dots, x_k) \rightarrow A(x_1, \dots, x_k) \, {\mod} \, p \Bigl \vert_{x_1=1, \dots, x_k=1}
\eqno(Reduce)
$$
to mean ``expand $A(x_1, \dots, x_k)$ as a sum of monomials,
reduce the coefficients modulo $p$ to one of the numbers  
$\{0, 1, \ldots, p-1\} \in {\Z}$,  and finally set all the variables $x_i$ equal to $1$''.

For any polynomial $Q=Q(x_1, \dots, x_k) \in {\Z}[x_1,\ldots,x_k]$ 
whose degree in each of the variables is less than $p$, define
$$
a_Q(n) := Q P^n \, {\mod} \, p \Bigl \vert_{x_1=1, \dots , x_k=1} \quad  .
$$
For $0 \leq i <p$, we have
$$
a_Q(pn\, + \, i) \, = \, Q(x_1, \dots, x_k)  P(x_1, \dots, x_k) ^{pn+i} \, {\mod} \, p \, \Bigl \vert_{x_1=1, \dots , x_k=1}
$$
$$
= \, [ \, Q(x_1, \dots, x_k)  P(x_1, \dots, x_k)^i \, ] P(x_1, \dots, x_k)^{np} \, {\mod} \, p \, \Bigl \vert_{x_1=1, \dots , x_k=1}
$$
$$
= \, [ \, Q(x_1, \dots, x_k)  P(x_1, \dots, x_k)^i \, ] (P(x_1, \dots, x_k)^{p})^{n} \, {\mod} \, p \, \Bigl \vert_{x_1=1, \dots , x_k=1}
$$
$$
= \, [ \, Q(x_1, \dots, x_k)  P(x_1, \dots, x_k)^i \, ] P(x_1^p, \dots, x_k^p)^{n} \, {\mod} \, p \, \Bigl \vert_{x_1=1, \dots , x_k=1} \quad .
$$

Now write
$$
 Q(x_1, \dots, x_k)  P(x_1, \dots, x_k)^i \, {\mod} \, p \, = \,
\sum_{(\alpha_1, \dots, \alpha_k) \in \{0,\dots , p-1\}^k } x_1^{\alpha_1} \cdots x_k^{\alpha_k} R_{(\alpha_1, \dots, \alpha_k)} (x_1^p, \dots, x_k^p)
\quad .
$$
(Here again ``mod $p$'' applies just to the coefficients, not the variables.)
Hence
$$
a_Q(np+i)=
\sum_{(\alpha_1, \dots, \alpha_k) \in \{0,\dots , p-1\}^k}
x_1^{\alpha_1} \cdots x_k^{\alpha_k} R_{(\alpha_1, \dots, \alpha_k)} (x_1^p, \dots, x_k^p) P(x_1^p, \dots, x_k^p)^n \, {\mod} \, p \, \Bigl 
\vert_{x_1=1, \dots, x_k=1}
$$
$$
=\, \sum_{(\alpha_1, \dots, \alpha_k) \in \{0,\dots , p-1\}^k } R_{(\alpha_1, \dots, \alpha_k)} (x_1^p, \dots, x_k^p) P(x_1^p, \dots, x_k^p)^n \, {\mod} \, p \, 
\Bigl \vert_{x_1=1, \dots , x_k=1} \quad,
$$
$$
=\, \sum_{(\alpha_1, \dots, \alpha_k) \in \{0,\dots , p-1\}^k } R_{(\alpha_1, \dots, \alpha_k)} (x_1, \dots, x_k) P(x_1, \dots, x_k)^n \, {\mod}\, p \,\Bigl 
\vert_{x_1=1, \dots , x_k=1} \quad,
$$
$$
=\, \sum_{(\alpha_1, \dots, \alpha_k) \in \{0,\dots , p-1\}^k } 
a_{R_{(\alpha_1, \dots, \alpha_k)}}(n) \quad .
$$
In other words for any $Q(x_1, \dots, x_k)$ and each of the residue classes 
$i$, $0\leq i \leq p-1$,
we can find a multiset of polynomials, let's call it $S_i(Q)$, such that
$$
a_Q(np+i) \, = \sum_{R \in S_i(Q)} a_R(n) \quad.
$$
We really only care about the case $Q=1$, but the algebra
 forces us to consider other $Q$'s, and they in turn force us to treat
still other $Q$'s, and so on. However, by the {\bf pigeon-hole principle},
this process must terminate, and we
obtain a {\bf finite} recurrence scheme, containing say $m$ equations.
Placing all the $Q$'s that appear into some arbitrary order,
with $Q_1=1$, we get a (logarithmic-time) {\bf recurrence scheme}:
$$
a_j(np+i) \, =\sum_{l \in S_i(j)} a_l(n) \quad ,
$$
for $1 \leq j \leq m$, that enables the fast calculation of $a_1(n)$ for any $n$.

Furthermore, by focusing only on $i=p-1$, 
and defining $c_j(k):=a_j(p^k-1)$, we have, for $1 \leq j \leq m$,
$$
c_j(k)\, =\sum_{l \in S_{p-1}(j)} c_l(k-1) \quad .
$$
Define the {\bf generating functions}
$$
f_j(t):=\sum_{k=0}^{\infty} c_j(k) t^k  \quad (1 \leq j \leq m) \quad .
$$
Standard manipulations of generating functions 
convert the above recurrences into a system of $m$ {\it linear} equations for the 
$m$ unknowns $f_1(t), \dots, f_m(t)$:
$$
f_j(t)= c_j(0)+ t \sum_{l \in S_{p-1}(j)} f_l(t) \quad , \quad 1 \leq j \leq m \quad ,
$$
that can be  solved, at least in principle, 
yielding {\bf rigorous} explicit expressions for {\it all} the $f_j(t)$,
and  in particular for  $f_1(t)$, the one in which we are most interested. 
Note that this proves that the generating function, $f_1(t)$, is always a {\bf rational function}.
If $m$ is too large, and the system of equations cannot be solved, then one may 
try to use the recurrences to  generate sufficiently many terms
of the sequence $c_1(k)$, and then
{\it guess} the rational function $f_1(t)$, 
using for example the Maple packgage {\tt gfun} [SaZ]. 
It may then be possible to 
justify that guess, {\it a posteriori}, by finding upper bounds on the degree
of the generating function.

{\bf Keeping track of the individual coefficients}

If instead of the functional Eq. ($Reduce)$, 
one uses, for some formal variables $s_1, \dots , s_{p-1}$, 
$$
\sum_{{\bf \alpha}} c_{{\bf \alpha}} {\bf x}^{{\bf \alpha}} \rightarrow \sum_{{\bf \alpha}} s_{c_{{\bf \alpha}}} \quad,
$$
one can modify the above arguments and 
keep track of the number of occurrences of each $i$ ($i=1, \ldots, p-1$)
as coefficients in the expansion of $P(x_1, \dots, x_k)^n \, {\mod} \, p$.

{\bf The Maple package CAcount}

Everything discussed above  is implemented in the Maple package {\tt CAcount}, which
can be down-loaded  from the web page for this article: \hfill\break
{\tt http://www.math.rutgers.edu/\~{}zeilberg/mamarim/mamarimhtml/CAcount.html} \qquad ,
where there are also  many samples of input and output files 
that readers can use as templates for further computations.

To see the list of the main procedures, type 

{\tt ezra();} \quad , 

or to see the
list of procedures that handle the more refined case, 
where one keeps track of the individual coefficients
(only useful for $p>2$), type 

{\tt ezraG();} \quad  .

To get instructions on using a particular procedure, type

{\tt ezra(ProcedureName);} \quad .

For example. procedure {\tt CAaut} finds the recurrence `automaton', and to get help with it,
type

{\tt ezra(CAaut);} \quad .

For our toy example, type

{\tt CAaut([1+x+x**2,1],[x],2,2);}

which produces as output the pair

{\tt [[[[1], [2, 1]], [[1, 1], [1, 1]]], [1, 2]]} \quad,

where the first component,

{\tt [[[1], [2, 1]], [[1, 1], [1, 1]]]} \quad ,

is Maple's way of encoding the recurrence
$$
a_1(2n)=a_1(n) \quad , \quad a_1(2n+1)=a_2(n)+a_1(n) \quad ;  \quad a_2(2n)=a_1(n)+a_1(n) \quad ,  \quad a_2(2n+1)=a_1(n)+a_1(n) \quad .
$$

The second component 

{\tt [1, 2]}

is Maple's way of encoding the initial conditions 

$$
a_1(1)= 1 \quad, \quad a_2(1)=2 \quad .
$$

Procedure {\tt SeqF} uses the scheme, once found, to compute as many terms as desired, while procedure {\tt ARLT}
(for {\it anti-run-length-transform}, see [Sl]) computes the sparse subsequence in the places $p^i-1$.
Procedure {\tt GFsP} finds the {\bf proved} generating function for that subsequence, and if the
size of the system is too big,  {\tt GFsG} guesses it faster, and as we mentioned above, the
guess can be justified {\it a posteriori}.

{\bf References}

[3by3] Shalosh B. Ekhad, N.~J.~A.~Sloane, and Doron Zeilberger,
{\it ``Odd-Rule'' Cellular Automata on the Square Grid}, in preparation, March 2015.

[OEIS] 
The OEIS Foundation Inc.,
{\it The On-Line Encyclopedia of Integer Sequences}, 
{\tt https://oeis.org}.

[RZ] Eric Rowland and Doron Zeilberger,
{\it A Case Study in Meta-AUTOMATION: AUTOMATIC Generation of Congruence AUTOMATA For Combinatorial Sequences},
J. Difference Equations and Applications {\bf 20} (2014), 973--988;
 \hfill\break
{\tt http://www.math.rutgers.edu/\~{}zeilberg/mamarim/mamarimhtml/meta.html}.

[SaZ] Bruno Salvy and Paul Zimmermann, {\it GFUN: a Maple package
for the Manipulation of Generating and Holonomic Functions
in One Variable}, ACM Trans. Math. Software {\bf 20} (1994), 163--177.

[Sl] N. J. A. Sloane {\it  On the Number of ON Cells in Cellular Automata}, 2015;
\hfill\break
{\tt http://arxiv.org/abs/1503.01168}.

\bigskip
\bigskip

\hrule
\bigskip
Shalosh B. Ekhad, c/o D. Zeilberger, Department of Mathematics, Rutgers University (New Brunswick), Hill Center-Busch Campus, 110 Frelinghuysen
Rd., Piscataway, NJ 08854-8019, USA.
\bigskip
\hrule
\bigskip
N. J. A. Sloane, The OEIS Foundation Inc, 11 South Adelaide Ave, Highland Park, NJ 08904, USA,  and Department of Mathematics, Rutgers University (New Brunswick);  \hfill \break
 njasloane at gmail dot com \quad ;  \quad {\tt  http://neilsloane.com/ } \quad .
\bigskip
\hrule
\bigskip
Doron Zeilberger, Department of Mathematics, Rutgers University (New Brunswick), Hill Center-Busch Campus, 110 Frelinghuysen
Rd., Piscataway, NJ 08854-8019, USA; \hfill \break
zeilberg at math dot rutgers dot edu \quad ;  \quad {\tt http://www.math.rutgers.edu/\~{}zeilberg/} \quad .
\bigskip
\hrule
\bigskip
\bigskip
Published in The Personal Journal of Shalosh B. Ekhad and Doron Zeilberger  \hfill \break
({\tt http://www.math.rutgers.edu/\~{}zeilberg/pj.html}), 
N.~J.~A.~Sloane's home page \hfill \break ({\tt  http://neilsloane.com/}),
and {\tt arxiv.org} \quad .
\bigskip
\bigskip
\hrule
\bigskip

{\bf Mar 05, 2015}

\end